\documentclass[10.5pt]{amsart}
\usepackage{amsmath,amssymb,amscd,enumerate,amsfonts}

\newtheorem{theorem}{Theorem}[section]

\newtheorem{corollary}[theorem]{Corollary}
\newtheorem{remark}[theorem]{Remark}

\usepackage{amsmath,amssymb,amscd,enumerate,amsfonts}

\numberwithin{equation}{section}

\begin{document}

\leftline{ \scriptsize {\it  The first draft} }

\bigskip
\bigskip

\title[Some Fractional Functional Inequalities]{Some Fractional Functional Inequalities and Applications to some Minimization constrained Problems involving a local linearity}

\author{Hichem Hajaiej}

\address
{Hichem Hajaiej:
 King Saud University, P.O. Box 2455, 11451 Riyadh, Saudi Arabia}
\email{\rm hichem.hajaiej@gmail.com}

\date{}

\begin{abstract}
....
\end{abstract}

\subjclass[2010]{....}

\maketitle

\section{Introduction}

 \vspace{1in}

 - introduction will be added -

 \vspace{1in}


The fractional Laplacian is characterized as
\[
\sqrt{-\Delta}^{\, s} \phi
:=  {\mathcal F}^{-1} (| \cdot |^s
{\mathcal F}(\phi)),
\]
where
$\hat{u}= {\mathcal F}(u)$ represents
the Fourier transform of $u$ on ${\mathbb R}^{n}$
defined by
\[
\hat{f}(\xi) = {\mathcal F}(f)(\xi)
=  \int_{{\mathbb R}^{n}} f(x) e^{- i x \cdot \xi}  \, dx,
\]
if $f \in L^1({\mathbb R}^{n}) \cap L^2({\mathbb R}^{n})$.

\section{Fractional integral inequalities and compact embedding}
In this section, we will construct the fractional Polya-Szeg\"o inequality,
and present a fractional version of Gargliardo-Nirenberg inequality.
As an application, we show that the fractional Sobolev space
$W^{s, p}(\mathbb{R}^n)$ is compactly embedded into Lebesgue spaces $L^q(\Omega)$.

\subsection{Fractional Polya-Szeg\"o inequality}
We investigate
the nonexpansivity of Schwarz symmetric decreasing rearrangement of functions
with respect to the fractional actions $(-\Delta)^{ s/2}$ for $s \geq 0$.
For the basic terminology and some properties of Schwarz
symmetric decreasing rearrangement,
we refer Chapter 3 in \cite{lieb-loss}, also \cite{Burchard-Hajaiej}.

\begin{theorem}
\label{FPS}
Let $0 \leq s \leq 1$.
Let $u^{*}$ denote the Schwarz symmetric radial decreasing rearrangement
of $u$. Then we have
\begin{align}
\int_{\mathbb{R}^n}
|  \sqrt{-\Delta}^{\, s} u^*(x)|^2 dx
\leq
\int_{\mathbb{R}^n}
|  \sqrt{-\Delta}^{\, s} u(x)|^2 dx,
\label{rearr_ineq}
\end{align}
in the sense that the finiteness of of the right side implies the finiteness
of the left side.
\end{theorem}

\proof
When $s = 0$, we have the equality in (\ref{rearr_ineq}).
We now are going to present how the kinetic energy decreases via
the symmetric radial decreasing rearrangement
as the differential index $s$ increases.

To show (\ref{rearr_ineq}),
it is enough to prove the following:
\begin{align}
 \int_{{\mathbb{R}}^{n}} |\xi|^{2s} |{\mathcal F}[f^{*}](\xi)|^{2} d \xi
\leq
 \int_{{\mathbb{R}}^{n}} |\xi|^{2s} |\hat f(\xi)|^{2} d \xi.
\label{FPS-eq*}
\end{align}
The main idea of the proof is that the inequality (\ref{FPS-eq*}) can be
followed from proving the assertion:
 for any $\varepsilon > 0$,
\begin{align}
 \int_{{\mathbb{R}}^{n}}
  \left( \frac{|\eta|^{2}}{1 + \varepsilon^{2}|\eta|^{2}}
  \right)^{s} |{\mathcal F}[u^{*}](\eta)|^{2} d \eta
\leq
 \int_{{\mathbb{R}}^{n}}
  \left( \frac{|\eta|^{2}}{1 + \varepsilon^{2}|\eta|^{2}}
  \right)^{s} |\hat u(\eta)|^{2} d \eta.
\label{equi*}
\end{align}
With change of variables $\xi = \varepsilon \eta$,
(\ref{equi*}) becomes
\begin{align}
 \frac{1}{\varepsilon^{2n}}\int_{{\mathbb{R}}^{n}}
    \left(\frac{|\xi|^{2}}{1 + |\xi|^{2}}
    \right)^{s}
     |{\mathcal F}[u^{*}](\xi/\varepsilon)|^{2} d \xi
\leq
 \frac{1}{\varepsilon^{2n}}\int_{{\mathbb{R}}^{n}}
    \left(\frac{|\xi|^{2}}{1 + |\xi|^{2}}
    \right)^{s}
     |\hat u(\xi/\varepsilon)|^{2} d \xi.
\label{Equ*}
\end{align}
Replace $u(x)$ by
$u(x / \varepsilon)$, and we have
$[u(x/ \varepsilon)]^{*} = u^{*}(x/ \varepsilon)$
since rearrangement commutes with uniform
dilation on the space.
Then   (\ref{Equ*}) is equivalent to saying
\begin{align}
 \int_{{\mathbb{R}}^{n}}
    \left(\frac{|\xi|^{2}}{1 + |\xi|^{2}}
    \right)^{s} |{\mathcal F}[u^{*}](\xi)|^{2} d \xi
\leq
 \int_{{\mathbb{R}}^{n}}
    \left(\frac{|\xi|^{2}}{1 + |\xi|^{2}}
    \right)^{s} |\hat u(\xi)|^{2} d \xi.
\label{Equ**}
\end{align}
So it suffices to prove (\ref{Equ**}).
Incorporating
the following expression
\begin{align*}
       \left(
         \frac{|\xi|^{2}}{1 + |\xi|^{2}}
       \right)^{s}
 =   \left(1 - \frac{1}{1 + |\xi|^{2}}
       \right)^{s}
 =   1 - \sum_{k = 1}^{\infty} (-1)^{k+1}
                \left( \begin{array}{c}
                        s \\ k
                       \end{array}
                \right)
             \left(\frac{1}{1 + |\xi|^{2}}
       \right)^{k}
\end{align*}
with
$
{s\choose k} = \frac{s(s - 1) \cdots (s-(k-1)) }
                                          {k!}
$
into each side of
inequality (\ref{Equ**}) yields
\begin{align*}
&\int_{{\mathbb{R}}^{n}} |{\mathcal F}[u^{*}](\xi)|^{2} d\xi -
\sum_{k = 1}^{\infty} (-1)^{k+1} \left( \begin{array}{c}
                        s \\ k
                       \end{array}
                \right)
                       \int_{{\mathbb{R}}^{n}}
                       \frac{1}{(1+|\xi|^{2})^{k}}
                       |{\mathcal F}[u^{*}](\xi)|^{2} d\xi
\end{align*}
and
\begin{align*}
\int_{{\mathbb{R}}^{n}} |{\mathcal F}[u](\xi)|^{2} d\xi -
\sum_{k = 1}^{\infty} (-1)^{k+1} \left( \begin{array}{c}
                        s \\ k
                       \end{array}
                \right)
                       \int_{{\mathbb{R}}^{n}}
                       \frac{1}{(1+|\xi|^{2})^{k}}
                       |{\mathcal F}[u](\xi)|^{2} d\xi.
\end{align*}
Since
$(-1)^{k + 1} \left( \begin{array}{c}
                        s \\ k
                       \end{array}
                \right) > 0$ with $0 < s < 1$,
it remains to show that for each positive integer $k$
\[
 \int_{{\mathbb{R}}^{n}}
                       \frac{1}{(1+|\xi|^{2})^{k}}
                       |{\mathcal F}[u^{*}](\xi)|^{2} d\xi
\geq
 \int_{{\mathbb{R}}^{n}}
                       \frac{1}{(1+|\xi|^{2})^{k}}
                       |{\mathcal F}[u](\xi)|^{2} d\xi.
\]
We consider a Bessel kernel $G_{2k}$ of order $2k$:
$
(1 + |\xi|^{2})^{-k} = {\mathcal F}[G_{2k}](\xi).
$
Therefore with $\tilde{u}(x)= u(-x)$, we arrive at
\begin{align}
\int_{{\mathbb{R}}^{n}}
           \frac{1}{(1+|\xi|^{2})^{k}}
           |{\mathcal F}[u](\xi)|^{2} d\xi
&   =  \int_{{\mathbb{R}}^{n}}
           {\hat G_{2k}}(\xi)
           \overline{{\hat u}}(\xi)
            {\hat u}(\xi) d\xi \nonumber \\
&   =   (2 \pi)^n \int_{{\mathbb{R}}^{n}}
           G_{2k}(-x)(u*\tilde{u})(x) dx \nonumber \\
&   =   (2\pi)^n [G_{2k}(x)*(u*\tilde{u})(x)](0)
              \nonumber \\
&   =   (2\pi)^n \int_{{\mathbb{R}}^{n} \times {\mathbb{R}}^{n}}
         G_{2k}(y-z)
         \bar{u}(z) u(y) dydz \nonumber \\
& \leq  (2\pi)^n \int_{{\mathbb{R}}^{n} \times {\mathbb{R}}^{n}}
         G_{2k}(y-z)
         \bar{u}^{*}(z) u^{*}(y)dydz
           \label{Re} \\
&   =   \int_{{\mathbb{R}}^{n}}
           \frac{1}{(1+|\xi|^{2})^{k}}
            |{\mathcal F}[u^{*}](\xi)|^{2} d\xi. \nonumber
\end{align}
The {\it Symmetrization lemma}
in \cite{Beckner1, lieb-loss} yields the inequality (\ref{Re})
where is the only place that inequality occurs.
The proof is now completed.
\hfill$\Box$\par

\subsection{Fractional Gargliardo-Nirenberg Inequality}

Gargliardo-Nirenberg inequality for fractional Laplacian is presented, and
sharp form of the fractional Sobolev inequality is obtained as a corollary.
Throughout this paper, $C$ denotes various real positive constants which do
not depend on functions in discussion.

\begin{theorem}
\label{FGN}
Let
$m, q, \theta \in \mathbb{R} \setminus \{ 0 \}$
with $q \neq m \theta > 0$, 
$0 < s < n$, $1 < p <\frac{n}{s}$
and $1 < \frac{r}{q-m\theta}$. Then the inequality
\begin{align}
\int_{\mathbb{R}^n}
| u(x)|^q dx
\leq C
\left(
\int_{\mathbb{R}^n}
\left( \sqrt{-\Delta}^{\, s} u(x) \right)^p  dx
\right)^{\frac{m \theta}{p}}
\left(
\int_{\mathbb{R}^n}
| u(x)|^r  dx
\right)^{\frac{q- m \theta}{r}}
\label{FGN-ineq}
\end{align}
holds for the indices with the relation
\begin{align}
m \theta
\left( \frac{1}{p}
- \frac{s}{n}
\right)
+ \frac{q - m \theta}{r} = 1.
\label{index}
\end{align}
In particular,
when $m = q$,  we have a fractional version of  Gargliardo-Nirenberg inequality:
\begin{align}
\left(
\int_{\mathbb{R}^n}
| u(x)|^q dx
\right)^{\frac{1}{q}}
\leq C
\left(
\int_{\mathbb{R}^n}
| \sqrt{-\Delta}^{\, s} u(x)|^{p} dx
\right)^{\frac{\theta}{p}}
\left(
\int_{\mathbb{R}^n}
| u(x)|^r  dx
\right)^{\frac{1- \theta}{r}}
\label{FGN-ineq1}
\end{align}
for the indices with the relation
\begin{align}
\theta \left( \frac{1}{p}
- \frac{s}{n}
\right)
+ \frac{1-\theta}{r} = \frac{1}{q}.
\label{index1}
\end{align}
\end{theorem}
\proof
For convenience, we use the notation
$\| u \|_{L^t} := \left( \int_{\mathbb{R}^n} | u(x)|^t  dx \right)^{\frac{1}{t}}$, and
$L^t(\mathbb{R}^n)$ for {\it any} $t \in \mathbb{R}\setminus \{ 0 \}$.
First we point out that
by the standard
{\it dilation argument} the index relation (\ref{index}) is necessary.
In fact,
by replacing
$u(\cdot)$ with $u(\delta  \; \cdot)$, we can observe
\[
\delta^{-n}
\| u \|_{L^q}^q
\leq C
\delta^{^{(s  - \frac{n}{p})\theta + \frac{n(m\theta -q)}{r}}  }
\| \sqrt{-\Delta}^{\, s} |u|^m \|_{L^p}^{\theta}
\| u \|_{L^r}^{q- m \theta},
\]
 for all $\delta>0$, which implies  that
 $-n = (s  - \frac{n}{p})\theta + \frac{n(m\theta -q)}{r}$.

Now,
for any $u \in {\mathcal S}(\mathbb{R}^n)$, we have
\begin{align}
\int_{\mathbb{R}^n} | u(x)|^q dx
&=
\int_{\mathbb{R}^n}
| u(x)|^{m \theta} | u(x)|^{q - m \theta} dx   \nonumber \\
&\leq
\|  |u|^{m \theta} \|_{L^{\bar{p}}}
\| |u|^{q- m \theta} \|_{L^{\bar{r}}},
           \qquad \frac{1}{\bar{p}} + \frac{1}{\bar{r}} = 1  \label{ineq1} \\
&=\| u \|_{L^{m \theta\bar{p}}}^{m \theta}
\| u \|_{L^{(q- m \theta)\bar{r}}}^{q- m \theta}.   \nonumber
\end{align}
We set $m \theta\bar{p} := p_0$ and $(q- m \theta)\bar{r} := r$ to have
\begin{align}
\int_{\mathbb{R}^n} | u(x)|^q dx
\leq
\| u \|_{L^{p_0}}^{m \theta}
\| u \|_{L^{r}}^{q- m \theta},
\label{eq1}
\end{align}
and
$\frac{m \theta}{p_0}
+ \frac{q - m \theta}{r} = 1$.
Let $\sqrt{-\Delta}^{\, s} u := f$, and we have
\begin{align*}
u(x)
&= \frac{c_{n-s}}{c_{s}}
\int_{\mathbb{R}^n} \frac{f(y)}{ |x - y|^{n - s } }  dy
= \frac{c_{n-s}}{c_{s}}
\int_{\mathbb{R}^n}
  \frac{\sqrt{-\Delta}^{\, s} u(y)}{ |x - y|^{n - s } }  dy,
\end{align*}
where
$
c_s = \frac{\Gamma(s/2)}{\pi^{s/2}}.
$
Indeed, we may  take the Fourier transform on
$\sqrt{-\Delta}^{\, s} u  = f$, and
take it back to have $u$ after solving for $\widehat{u}$.
Therefore the Hardy-Littlewood-Sobolev inequality yields
\begin{align}
\| u \|_{L^{p_0}}
&\leq
\frac{c_{n-s}}{c_{s}}
C_1
\left\|
          \sqrt{-\Delta}^{\, s } u
\right\|_{L^{p}},
\label{eq2}
\end{align}
where $C_1$ is a positive constant(see the remark after the proof) and
$p$ satisfies
\begin{align}
\frac{1}{p} + \frac{n - s}{n} = 1 + \frac{1}{p_0}.
\label{index3}
\end{align}
This index relation combining with the index relation appeared at  (\ref{eq1})
implies  (\ref{index}),
and (\ref{eq1}) together with (\ref{eq2}) implies (\ref{FGN-ineq}).
\hfill $\square$  \\

It is known the best constant $C_1$ and the extremals  of the
Hardy-Littlewood-Sobolev inequality
for some special cases(see \cite{Lieb1}
or Section 4.3 in \cite{lieb-loss}).
Thanks to those cases, we have a sharp form of the fractional Sobolev inequality:
\begin{corollary}[Fractional Sobolev inequality]
\label{2-coro1}
For
$0 < s < n$, $1 < p <\frac{n}{s}$
and $q = \frac{pn}{n - sp}$, we have
\[
\| u \|_{L^q}
\leq C_0
\left\| \sqrt{-\Delta}^{\, s}  u \right\|_{L^p}.
\]
The sharp constant for the inequality  is
\[
\pi^{s/2} \frac{\Gamma(\frac{n-s}{2})}{\Gamma(\frac{n+s}{2})}
\left\{
  \frac{\Gamma(n)}{\Gamma(\frac{n}{2})}
\right\}^{s/n}.
\]
\end{corollary}
For a special case, we emphasize  the $L^2$-estimate of
the fractional Gargliardo-Nirenberg inequality which is applied at Section \ref{FSeq}.
\begin{corollary}
\label{2-coro}
For $0 < s < \frac{n}{2}$, $0 < \theta <1$ and $\theta = \frac{n(q-2)}{2qs}$,
we have
\[
\| u \|_{L^q}
\leq C
\left\| \nabla_s  u \right\|_{L^2}^{\theta}
\left\| u \right\|_{L^2}^{1- \theta}
\]
with the notation
$\| \nabla_s   u  \|_{L^2}
  :=
                            \left(
                             \int_{\mathbb{R}^n} |
                                       (- \Delta)^{\frac{s}{2}} u(x)
                                    |^2 dx
                             \right)^{\frac{1}{2} }$.
\end{corollary}

\subsection{Fractional Rellich-Kondrachov Compactness theorem}

The following theorem illustrates that the fractional Sobolev space
$W^{s, p}(\mathbb{R}^n)$ is compactly embedded into Lebesgue spaces $L^q(\Omega)$,
where $\Omega$ is bounded.
\begin{theorem}
\label{R-C-thm}
Let $0 < s < n$,
$1 \leq  p <\frac{n}{s}$  and $1 \leq q < \frac{np}{n - sp}$.
Also, let $\{ u_m  \}$ be a sequence in $L^q(\mathbb{R}^n)$
and $\Omega$ be a bounded open set
with smooth boundary.
Suppose that
\begin{align*}
\int_{\mathbb{R}^n}
| \sqrt{-\Delta +1}^{\, s} u_m(x)|^{p} dx
\end{align*}
are uniformly bounded, then $\{ u_m \}$ has a convergent subsequence in $L^q(\Omega)$.
\end{theorem}
\proof
Let $\phi$ be a smooth non-negative function with support in
$\{ x  : |x| \leq 1 \}$ and with $\int_{|x| \leq 1}  \phi(x) dx  = 1$.
We also define $\phi^{\ell}(x) := {\ell}^n \phi(\ell x)$.
By virtue of the Fractional Sobolev inequality(Corollary \ref{2-coro1}),
it can be observed that
\begin{align}
\| u_m  \|_{L^q(\mathbb{R}^n)}
\leq C
\left\| \sqrt{-\Delta}^{\, s}  u_m \right\|_{L^p(\mathbb{R}^n)}
\leq C
\left\| \sqrt{-\Delta +1}^{\, s}  u_m \right\|_{L^p(\mathbb{R}^n)}
\leq \widetilde{C}    \label{2-est1}
\end{align}
for some $\widetilde{C} > 0$.
Hence in the spirit of Frechet-Kolmogorov theorem,
it suffices to show the following (see page 50 in \cite{show}):
for any $\varepsilon >0$ and
any compact subset $K$ of $\Omega$,
there is a constant $M >0$ such that
for $m \geq M$,
\[
\| \phi^{\ell} *u - u \|_{L^q(K)}  < \varepsilon,
\]
for all  $u \in \mathcal{S}(\mathbb{R}^n)$ with
$\| \sqrt{-\Delta +1}^{\, s} u \|_{L^p(\mathbb{R}^n)}
  \leq \widetilde{C}/ C$.
Then using the interpolation inequality (\ref{eq1}), we have
\begin{align*}
\| \phi^{\ell} *u - u \|_{L^q(K)}
\leq C 2^{1-\theta}
\| u \|_{L^r(\mathbb{R}^n)}^{1-\theta}
\| \phi^{\ell} *u - u \|_{L^1(K)}^{\theta},
\end{align*}
with $\frac{1- \theta}{r} + \theta = \frac{1}{q}$,
$r = \frac{np}{n - sp}$.
Consequently, (\ref{2-est1}) and the fractional
Gargliardo-Nirenberg inequality(Theorem \ref{FGN})
imply that
\begin{align*}
\| \phi^{\ell} *u - u \|_{L^q(K)}
\leq C
\| \phi^{\ell} *u - u \|_{L^1(K)}^{\theta}.
\end{align*}

Now we define $f := \sqrt{-\Delta+1}^{\, s} u$ to have
$u = G_s  * f$ and $\| f \|_{L^p} \leq \widetilde{C}/ C$,
where
$G_s$ is the Bessel kernel of order $s$.
Therefore we obtain
\[
\| \phi^{\ell} *u - u \|_{L^1(K)}
\leq C \| (\phi^{\ell} *G_s - G_s)*f \|_{L^p(\mathbb{R}^n)}
\leq C \| \phi^{\ell} *G_s - G_s \|_{L^1(\mathbb{R}^n)} \to 0
\]
as $m \to \infty$.
\hfill$\Box$\par

\section{Ground state solution of fractional Schr\"odinger flows}
\label{FSeq}

We consider the following variational problem:
\begin{align}
I_c
     &= \inf
        \left\{
           E(u) :  u\in S_c
        \right\}      \label{3-Pc}
\end{align}
where
$c$ is a prescribed number, $0<s<1$
and $E$ is the energy functional
\[
E(u) = \int_{{\mathbb{R}}^n}
                   | \sqrt{-\Delta}^{\, s} u(x)|^{2} dx
        - \int_{{\mathbb{R}}^n}
                  F( |x|, u(x) ) dx
\]
on an admissible collection
$
S_c :=\left\{
          u \in H^{s}(\mathbb{R}^n) :  \int_{\mathbb{R}^n} {u^{2}}(x) \, dx = c^2
     \right\}
$.

The aim of this work is to study the symmetry properties of minimizers of (\ref{3-Pc}).
We can also note that the solutions of (\ref{3-Pc}) lie on the curve
\begin{align}
 (-{\Delta})^{s}u + f(|x|, u ) +  \lambda u = 0,  \label{3-pde}
\end{align}
where
$\lambda$ is a Lagrange multiplier  and
$
F(r, s)= \int_{0}^{s} f(r, t)dt.
$
It will be 
interesting to study the above identity (\ref{3-pde}) and to find
suitable assumptions on $f$ for which all the solutions of (\ref{3-pde})
are radial and radially decreasing.
Note that for the classical Laplacian, H. JeanJean(?) and C. Stuart have
completely  solved the problem.
It is also worth to study a fractional Schr\"odinger equation
\begin{align}
\left\{
 \begin{array}{cc}
        i\partial_{t} \Phi + (- {\Delta})^{s} \Phi + f(|x|, \Phi ) = 0     \\
        \Phi (x,0) = \Phi_0 (x)
 \end{array}
\right.    \label{3-E}
\end{align}
for which ground state solutions $u$ of (\ref{3-pde}) give rise
to ground state solitary wave $\Phi$ of (\ref{3-E}).
The minimizing problem (\ref{3-Pc}) that we are going to look at imposes the following
assumptions on the function $F$:    \\
$(F_0)$
$F:[0,\infty) \times \mathbb{R} \rightarrow \mathbb{R}$
is a Carath\'eodory function, that is to say: \\
\indent \hspace{.1in} $\bullet $
                   $F(\cdot , s): [0, \infty) \rightarrow \mathbb{R}$
is measurable for all $s \in \mathbb{R}$  and \\
\indent \hspace{.1in} $\bullet $
$F(r ,\cdot ): \mathbb{R} \rightarrow \mathbb{R}$
is continuous for almost every $r \in [0,\infty)$.   \\
$(F_1)$
$F(r , s)   \leq    F(r ,| s | ) $
for almost every  $r  \geq  0 $ and all $s \in \mathbb{R}$.    \\
$(F_2)$
There are $ K > 0 $ and $ 0 < l < \frac{4s}{n} $ satisfying
for any $r, s \geq 0$,
\[
0 \leq   F(r , s)  \leq  K( s^2 + s^{l+2}).
\]
$(F_3)$
For every $ \varepsilon > 0 $, there exist $R_0, s_0 > 0$ such that
$F(r , s) \leq \varepsilon |s|^2 $
for almost every $r  \geq R_0$ and all $ 0 \leq s < s_0 $.    \\
$(F_4)$
The mapping $ (t,y)\mapsto F(  \frac{1}{t} , y ) $
is super-modular on $\mathbb{R}_{+} \!\!\! \times \mathbb{R}_{+} $, in other words,
\[
F(r , a ) + F(R ,A) \geq F( r, A) + F(R ,a)
\]
for all $r < R$ and $ a < A$.

\begin{theorem}
Under the conditions $(F0) \sim (F4)$,
the minimizing problem (\ref{3-Pc})   admits
 a Schwarz symmetric minimizer for any fixed constant $c$.
  Moreover if $(F4)$ holds with a
strict sign, then for any c, all minimizers of (\ref{3-Pc})  are Schwarz symmetric.
\end{theorem}

 A Schwarz symmetric function is a radial decreasing function.
 For more detailed accounts, we refer [BH].

\proof
1. {\it  Well-posedness of the problem (\ref{3-Pc})} (that is, $I_c > - \infty$):
We first show that
all minimizing sequences are bounded in $H^s(\mathbb{R}^n)$.
By $(F1)$ and $(F2)$, we can write
\begin{align*}
\int F(|x|, u(x))dx  \leq  \int F \left(|x|, |u(x)|  \right)dx
                                 \leq  K c^2  + K
                                                 \int |u(x)|^{l+2} dx.
\end{align*}
By virtue of the fractional Gagliardo-Nirenberg inequality(Corollary \ref{2-coro})
and Young's inequality,
there exists constant $K'$ such that
\begin{align}
\int_{\mathbb{R}^n} |u(x) |^{l+2} dx
&\leq
                                K' \biggl\{ \int_{\mathbb{R}^n}
                                       u^2(x) dx
                                      \biggl\} ^{(1-\theta)\frac{(l+2)}{2}}
                                       \| \nabla_s   u \|_{L^2}^{\theta(l+2)}.
                                       \label{3-1-1}   \\
&\leq  K' \frac{\varepsilon^p}{p}
                                   \biggl\{
                                       \| \nabla_s   u \|_{L^2}^2
                                   \biggl\}^{p \theta \frac{(l+2)}{2}}
       \!\! +
                               \frac{K'}{q \varepsilon^{q}}
                                   \biggl\{
                                          \int_{\mathbb{R}^n} u^2(x) dx
                                   \biggl\}^{q (1-\theta)\frac{(l+2)}{2}}  \nonumber
\end{align}
for any $\varepsilon > 0$, $p >1 $,  $\frac{1}{p}+\frac{1}{q} = 1$ and
$
\theta = \frac{nl}{2s(l+2)}
$.
We choose
$p = \frac{2}{\theta (l+2)} = \frac{4s}{nl} $ to get
\begin{align*}
\int_{\mathbb{R}^n} |u(x) |^{l+2} dx & \leq  \frac{K'}{p} \varepsilon^p
                                             \biggl\{
                                                     \| \nabla_s   u \|_{L^2}^2
                                             \biggl\}
       + \frac{K'}{q \varepsilon^{q} }
                                              \biggl\{
                                                     \int_{\mathbb{R}^n} u^2(x) dx
                                               \biggl\}
                                                     ^{q(1-\theta)\frac{l+2}{2}}   \\
                                     & = \frac{K'}{p} \varepsilon^p
                                           \| \nabla_s   u \|_{L^2}^2
                                        +  \frac{K'}{q \varepsilon^{q} }
                                           c^{q(1-\theta)(l+2)}.
\end{align*}
Therefore applying $(F2)$, we conclude
\begin{align*}
E(u) & \geq \frac{1}{2} \| \nabla_s   u \|_{L^2}^2
           - K c^2
           - \frac{K' K}{p} \varepsilon^p \| \nabla_s   u \|_{L^2}^2
           - \frac{K' K }{q \varepsilon^q} c^{q(1-\theta)(l+2)}    \\
     & = \left( \frac{1}{2} - \frac{K' k }{p} \varepsilon^p \right)
                       \| \nabla_s   u \|_{L^2}^2
        - K c^2
        -  \frac{K' K }{q \varepsilon^q} c^{q(1-\theta)(l+2)}.
\end{align*}

\begin{remark}
1. If we allow $l = \frac{4s}{n} $ in $(F2)$, the problem (\ref{3-Pc})  still makes sense for
sufficiently small values of $c$. In fact,
with $\theta = \frac{2}{l+2}$  and in view of (\ref{3-1-1})
we have
\[
\int_{\mathbb{R}^n} |u(x) |^{l+2} dx
   \leq  K' c^{\frac{4s}{n}}
         \| \nabla_s   u \|_{L^2}^2
\]
for $u \in S_c$.
Hence we get
\begin{align*}
E(u) & \geq  \frac{1}{2}
                      \| \nabla_s   u \|_{L^2}^2
             - K c^2
             - K' K c^{\frac{4s}{n}}  \| \nabla_s   u \|_{L^2}^2             \\
     & = \left(
           \frac{1}{2} - K' K c^{\frac{4s}{n}}
          \right)
          \| \nabla_s   u \|_{L^2}^2
        - K c^2.
\end{align*}
Thus $I_c > - \infty $ and all minimizing sequences are bounded in $H^{s}(\mathbb{R}^n)$
 provided that $0 < c < (\frac{1}{2K K'})^{\frac{n}{4s}}$.

2. We can prove that $I_c = - \infty$ for $l > \frac{4s}{n}$.
\end{remark}

2. {\it Existence of a Schwarz symmetric minimizing sequence.}
First note that if $u \in H^{s}(\mathbb{R}^n) $,    then
$ |u| \in H^{s}(\mathbb{R}^n)$.
In view of $(F1)$,  we certainly have that
\[
E(|u|) \leq E(u),  \quad  \mbox{for all } u \in  H^{s}(\mathbb{R}^n).
\]
Now by virtue of the fractional Polya-Szeg\"o inequality(Theorem \ref{FPS}):
\[
\| \nabla_s  |u|^{*} \|_{L^2}
\leq
\|  \nabla_s   |u| \|_{L^2}
\]
and Theorem 1 of $[BH]$,  we can observe that
\[
\int_{\mathbb{R}^n} F (
         |x|, |u|(x)
        )  dx
\leq
\int_{\mathbb{R}^n}  F (
         |x|, |u|^{*}(x)
        )  dx.
\]
Thus, without loss of generality, we may say that
(\ref{3-Pc})  always admits a Schwarz symmetric minimizing sequence.

3. {\it Let $\{ u_m \} = \{ u_m^* \} $  be a  Schwarz symmetric minimizing sequence. If
$\{ u_m \}$ converges weakly to  $u$ in $H^s(\mathbb{R}^n)$, then
\[
E(u) \,\, \leq \,\, \liminf_{m \to \infty} E\{ u_m \}.
\]}
The weak lower semi-continuity of  $L^2$-norm yields
\[
\| \nabla_s u \|_{L^2} \,\,
                                    \leq  \,\,
\liminf_{m \to \infty} \| \nabla_s u_m \|_{L^2}.
\]
Hence the assertion will follow by showing that
\[
\lim_{m\rightarrow \infty}
\int_{\mathbb{R}^n} F(|x|, u_m(x) ) dx
=
\int_{\mathbb{R}^n} F(  |x|, u(x) )  dx.
\]
For $R > 0$, let us first prove that:
\[
\lim_{m \rightarrow \infty}
                  \int_{|x| \leq R }
                      F(|x|, u_m (x) ) dx
=
\int_{|x| \leq R }
  F( |x|, u(x) ) \, dx.
\]
By the fractional Rellich-Kondrachov theorem(Theorem \ref{R-C-thm}),
$\{ u_m \}$ converges strongly to u in $L^{l+2}(\{x : |x| \leq R \})$.
Thus there exists a subsequence $\{u_{m_k}\}$ of $\{ u_m \}$ such that
$u_{m_k}(x) \rightarrow u(x)$  for almost every  $|x| \leq R$
and there is $h \in L^{l+2}(\{x : |x| \leq R \})$ satisfying
$|u_{m_k}| \leq h $.
We apply $(F2)$ to have
\[
F(|x|, u_{m_k}(x)) \leq K ( h^2(x) + h^{l+2}(x) ).
\]
Noticing that $h^2 + h^{l+2} \in L^{1}(\{x : |x| \leq R \})$,
the dominated convergence theorem gives
\[
\lim_{m \rightarrow \infty}
                  \int_{|x| \leq R }
                      F(|x|, u_m(x) ) dx
=
\int_{|x| \leq R }
  F(|x|, u(x) )  \, dx.
\]
Since $ u_m = u_m^*$, we now have
\[
\omega_n |x|^{n} u^2_m(x)
                                              \leq
\int_{|y| \leq |x| } u^2_m(y) dy
                                              \leq
c^2,
\]
where $\omega_n$ is the measure of the $n$-dimensional unit ball.
Thus we get
\[
u_m(x) \leq
         \frac{c}{ { \omega_n^{\frac{n}{2}} }  |x|^{ \frac{n}{2} }}
       \leq
         \frac{c}{ { \omega_n^{\frac{n}{2}} }  R^{ \frac{n}{2} }},
\,\,\,\,\,\,\, \mbox{for all } |x| > R.
\]
Therefore
for $ \varepsilon > 0$ and $R$ sufficiently large, we obtain by using $(F3)$ that
\[
\int_{|x| > R }
      F( |x|, u_m(x) ) dx
                                  \leq
\varepsilon \int_{|x| > R }
     u^2_m(x)  dx
                                  < \varepsilon c^2,
\]
which in turn implies that
$
\lim_{R \rightarrow \infty}
\lim_{n\rightarrow \infty}
                  \int_{|x| > R }
                      F(|x|, u_m(x) ) dx =0.
$
Since $u$ inherits all the properties used to get the above limit,
 it follows also that
\[
\lim_{R \rightarrow \infty}
\int_{|x| > R }
  F( |x|, u(x) ) \, dx = 0.
\]

4. We claim that $ u \in S_c$. Notice that
$S_c = H^s(\mathbb{R}^n) \cap \Lambda^{-1}(\{ c \})$, where
$\Lambda$ is
defined by $\Lambda(u) := \| u \|_{L^2}$ for $u \in L^2(\mathbb{R}^n)$.
We choose a Schwarz symmetric minimizing sequence $\{ u_m \} \subset S_c$
converging weakly to $u$ in $H^s(\mathbb{R}^n)$,
and so it converges strongly to $u$ in $L^2(\mathbb{R}^n)$.
Hence we have that  $u \in H^s(\mathbb{R}^n)$ and
$u \in \Lambda^{-1}(\{ c \})$.
Indeed, since $\Lambda$ is continuous,
$\Lambda^{-1}(\{ c \})$ is a closed set in $L^2(\mathbb{R}^n)$.

\vspace{.5in}

 - References should be added and replaced -


\begin{thebibliography}{10}



\bibitem{Beckner*}
W. Beckner,
\emph{Inequalities in Fourier analysis, }
Ann. Math., 102 : 159-182, 1975.

\bibitem{Beckner**}
W. Beckner,
\emph{Geometric inequalities in Fourier analysis,
Essays on Fourier Analysis in honor of Elias M. Stein, }
Princeton University Press, 36-68, 1995.

\bibitem{Beckner1}
 W. Beckner,
\emph{Sobolev inequalities, the Poisson semigroup, and
analysis on the sphere on ${\bf S}^{n}$, }
Proc. Natl. Acad. Sci., 89 : 4816-4819, 1992.

\bibitem{Burchard-Hajaiej}
A. Burchard, H. Hajaiej,
\emph{Rearrangement inequalities for functional with monotone integrands},
J. Funct. Anal. 233, 561-582, 2006.


\bibitem{Hajaiej}
H. Hajaiej and C.A. Stuart,
\emph{Existence and non-existence of Schwarz symmetric ground states for
elliptic eigenvalue problems, }
Ann. Mat. Pura Appl., 184 : 297-314, 2005

\bibitem{Kawohl}
B. Kawohl,
\emph{Rearrangements and convexity of level sets in PDE, }
Springer-Verlag, 1985.

\bibitem{Lieb}
E.H. Lieb,
\emph{Existence and uniqueness of the minimizing solution of
Choquard¡¯s nonlinear equation, }
Studies in Appl. Math. 57 : 93-105, 1976/77.

\bibitem{Lieb1}
E. H. Lieb,
\emph{Sharp constants in the Hardy-Littlewood-Sobolev and related inequalities},
Ann. of Math. 118 : 349-374, 1983.

\bibitem{lieb-loss}
E.H. Lieb and M. Loss,
\emph{Analysis, }
Volume 14 of Graduate Studies in Mathematics, AMS, 1997.

\bibitem{Polya-Szego}
G. Polya and G. Szeg\"o,
\emph{Inequalities for the capacity of a condenser, }
Amer. J. Math. 67: 1-32, 1945.

\bibitem{show}
R. E. Showalter,
\emph{Monotone Operators in Banach Space
and Nonlinear Partial Differential Equations},
volume 49 of Math. Surveys and Monographs, AMS, 1997.

\end{thebibliography}
\end{document}